\documentstyle[10pt,twoside]{siamltex}
\date{}

\newtheorem{teo}[theorem]{Theorem}
\newtheorem{lem}[theorem]{Lemma}
\newtheorem{cor}[theorem]{Corollary}
\newtheorem{pro}[theorem]{Proposition}
\newtheorem{rem}[theorem]{Remark}

\newtheorem{num}[theorem]{}
\def\dem{{\it Proof.\ }\rm}

\newfont{\bb}{msbm10}
\def\Bbb#1{\mbox{\bb #1}}

\def\Bbb#1{\mbox{\bb #1}}

\def\zR{\Bbb R}

\def\TT{{\rm T} \hskip -5pt{\rm T} \hskip4pt}
\def\bull{\vrule height 1.0ex width .4ex depth -.1ex }
\def\inc{\subseteq}
\def\QED{\bull}
\def\inv{^{-1}}
\def\*a{\#\sb a}

\def\la{\lambda}
\def\eps{\varepsilon}
\def\fii{\varphi }
\def\fia{\varphi _a}
\def\fib{\varphi _b}
\def\uni{{\cal U}}
\def\unia{{\cal U}_a}
\def\uniA{{\cal U}_A}
\def\unib{{\cal U}_b}
\def\Gs{G^s}
\def\Gpo{G^+}
\def\ca{{\cal A}}
\def\cH{{\cal H}}

\def\rai{^{1/2}}
\def\mrai{^{-1/2}}
\def\pim{\pi^+}
\def\piA{\pi_A}
\def\api{\langle}
\def\cpi{\rangle}
\def\pou{\Theta}
\def\noi{\noindent}
\def\csta{C$^*$-algebra}
\def\noi{\noindent}

\def\csta{C$^*$-algebra}

\def\ben{\begin{enumerate}}
\def\een{\end{enumerate}}
\def\beq{\begin {equation}}
\def\eeq{\end {equation}}
\def\barr{\begin{array}}
\def\earr{\end{array}}
\begin{document}
\title{ Polar decomposition under perturbations of the scalar product{}}
\author {G. Corach \thanks {Instituto Argentino de Matem\'atica, 
Saavedra 15 Piso 3 (1083), Buenos Aires Argentina 
(e-mail: gcorach@mate.dm.uba.ar). Partially supported by Fundaci\'on Antorchas, 
CONICET (PIP 4463/96), Universidad de Buenos Aires (UBACYT TX92 and TW49) and 
ANPCYT (PICT97-2259)} \and  A. Maestripieri \thanks{ Instituto de Ciencias, UNGS, 
Roca 851, San Miguel, Buenos Aires, Argentina (e-mail: amaestri@ungs.edu.ar)}
\and D. Stojanoff
\thanks {Depto. de Matem\'atica, UNLP, 1 y 50 (1900), La Plata, Argentina  
(e-mail: demetrio@mate.dm.uba.ar). Partially supported CONICET (PIP 4463/96), 
Universidad de Buenos Aires (UBACYT TW49) and ANPCYT (PICT97-2259)}}

\maketitle 
\vskip.5truecm
\centerline{Dedicated to our friend and teacher  Angel Rafael Larotonda}
\centerline{in his sixtieth anniversary}
\vskip.5truecm

\begin{abstract}{
Let ${\cal A}$ be a unital \csta \ with involution $*$ represented 
in a Hilbert space $\cH$, 
$G$ the group of invertible elements  of $\ca$, 
$\uni$ the unitary group of ${\cal A}$, $\Gs$ 
the set of invertible selfadjoint elements of $\ca$, 
$Q= \{ \eps \in G : \eps ^2 = 1 \}$ the space of reflections
and $P = Q \cap \uni$.  For any positive
$a \in  G$ consider the $a$-unitary group 
$\unia = \{ g \in G : a\inv g^* a = g \inv \}$, i.e. the elements 
which are unitary with respect to the scalar product 
$\langle \xi , \eta \rangle \sb a = \langle a\xi , \eta \rangle $ for 
$ \xi , \ \eta \in \cH $. If 
$\pi $ denotes the map  that assigns to each invertible
element its unitary part in the polar decomposition, 
we show that the restriction 
$\pi | _{\unia} : \unia \to \uni $
is a diffeomorphism, that $\pi ( \unia \cap Q) = P $ and 
that $ \pi (\unia \cap \Gs ) = \unia \cap \Gs = \{ u \in G : 
u = u^* = u \inv \  \hbox{ and} \ au = ua \} .$
}\end{abstract}

\begin{keywords}
Polar decomposition, \csta s, positive operators, projections.
\end{keywords}

\begin{AMS}
Primary 47A30, 47B15 \end{AMS}

\vskip 1truecm
\section{Introduction}
If $\ca $ is the algebra of bounded linear operators in a Hilbert space $\cH$, 
denote by $G$ the group of invertible elements of $\ca$. Every $T\in G$ admits two
{\bf polar decompositions \rm} 
$$
T= U_1 P_1 = P_2 U_2
$$
where $U_1, U_2$ are unitary operators (i.e. $U_i ^* = U_i \inv$) and
$P_1 , P_2 $ are positive operators (i.e $\api P_i \xi , \xi \cpi \ge 0$ for
every $\xi \in \cH$). It turns out that $U_1 = U_2$, $P_1 = (T^*T)\rai $ and
$P_2 = (TT^*)\rai$. We shall call $U = U_1 = U_2$ {\bf the unitary part \rm} 
of $T$. Consider the map
\beq \label{DP}
\pi : G \to \uni \quad \quad \pi (T) = U
\eeq
where $\uni$ is the unitary group of $\ca$. If $\Gpo $ denotes the set of all 
positive invertible elements of $\ca$, then every $A\in \Gpo$ defines an 
inner product $\api \ , \ \cpi _A$ on $\cH$ which is equivalent to the 
original $\api \ , \ \cpi $; namely
\beq
\api \xi , \eta \cpi _A = \api A \xi , \eta \cpi \quad (\xi , \eta \in \cH ).
\eeq
Every $X \in \ca$ admits an $A$-adjoint operator $X^{*_A}$, which is the unique 
$Y \in \ca$ such that 
$$
\api X \xi , \eta \cpi _A = \api \xi , Y \eta \cpi _A \quad (\xi , \eta \in \cH ).
$$
It is easy to see that $X^{*_A} = A\inv X^* A$. 
Together with the definition of $*_A$ one gets
the sets of {\bf $A$-Hermitian operators \rm} 
$$
\ca_A^h = \{ X \in \ca : X ^{*_A} = X  \} =\{ X \in \ca : AX  = X^*A \},
$$
{\bf $A$-unitary operators \rm} 
$$
\uniA =  \{ X \in G : X ^{*_A} = X\inv  \} =\{ X \in \ca : AX\inv  = X^*A \}
$$
and {\bf $A$-positive operators \rm}
$$
\Gpo_A = \{ X\in \ca_A^h \cap G : \api X \xi , \xi \cpi _A \ge 0 \quad
\forall \  \xi \in \cH \}
$$
As in the ``classical'' case, i.e. $A=I$, we get  polar decompositions
of any $T \in G$ 
$$
T=V_1R_1 = R_2 V_2 
$$
with $V_i \in \uniA $, $R_i \in \Gpo _A $, $i = 1,2$ and as before $V_1 = V_2$. 
Thus, we get a map
\beq
\piA : G \to \uniA .
\eeq
This paper is devoted to a simultaneous study of the maps 
$\piA$ ($A\in \Gpo$), the way that 
$$
\uniA \ , \quad G_A^h \ , \quad \Gpo _A 
$$
intersect
$$
\uni_B\  ,  \quad G_B^h \ , \quad \Gpo _B 
$$
for different $A, B \in \Gpo$ and the intersections of these
sets with
$$
Q = \{ S\in \ca : S^2 = I \}  
\quad \hbox { and } \quad
P_A = \{ S \in Q : S^{*_A} = S \}
$$
(reflections and  $A$-Hermitian reflections of $\ca$). 
The main result is the fact that, 
for every $A, B \in \Gpo$, 
$$
\pi _A |_{\uni _B} : \uni _B \to \uniA 
$$
is a bijection. The proof of this theorem is based on the form of 
the positive solutions of the operator equation
$$
XAX = B
$$
for $A, B \in \Gpo$. This identity  was first studied by G. K. Pedersen and 
M. Takesaki \cite{[PT1]} in their study of the 
Radon-Nykodym theorems in von Neumann algebras. As a corollary 
we get a short proof of the equality 
\beq
\pi_A (Q \cap \uni_B ) = P_A
\eeq
for every $A, B \in \Gpo$, which was proven in \cite{[ACS2]} 
as a $C^*$-algebraic version of results of 
Pasternak-Winiarski \cite{[PW1]} on the analyticity of the map 
$A \mapsto P_A^{\cal M}$, where 
$P_A^{\cal M}$ is the $A$-orthogonal projection on the closed subspace ${\cal M}$ 
of $\cH$. We include a parametrization of all solutions of Pedersen-Takesaki
equation.
The results are presented in the context of unital \csta s.

\section{Preliminaries}

Let ${\cal A}$ be a unital \csta ,  
$G= G(\ca )$ the group of invertible elements  of $\ca$, 
$\uni = \uni (\ca ) $ the unitary group of ${\cal A}$, $\Gpo = \Gpo (\ca )$  the 
set of positive invertible elements of $\ca$ and $\Gs \Gs (\ca )$ the set of positive 
selfadjoint elements of $\ca$. Let $Q = Q( \ca )= \{ \eps \in G : \eps ^2 = 1 \}$  the space
of reflections and 
$$
P = P( \ca ) = Q \cap \Gs = Q \cap \uni = \{ \rho \in G : \rho = \rho ^* = \rho \inv \}
$$
the space of orthogonal reflections, also called 
the Grassmann manifold of ${\cal A}$.

Each $g\in G$ admits two \bf polar decompositions \rm
$$
g = \la u = u'\la ' \ , \quad \la , \la'  \in \Gpo \ , \quad u, u'\in \uni .
$$
In fact, $\la = (gg^*)\rai $, $ u = (gg^*)\mrai g$, $\la'= (g^*g)\rai $ and 
$u'= g(g^*g)\mrai $. A simple exercise of functional calculus shows that $u = u '$.
We shall say that $u $ is $ the$  $unitary$  $part$ of $g$.
Observe that in the decomposition $g = \la u $ (resp. $g = u \la '$) the components
$\la , \ u$ (resp. $u$, $\la '$) are uniquely determined, for instance,
if $\la u = \la _o u _0 $, then $\la _0 \inv \la = u_0 u\inv $ is a 
unitary element wuth positive spectrum: $\sigma (\la _0 \inv \la ) = 
\sigma (\la _0 \mrai \la \la _0 \rai ) = \sigma (\la ) \inc \zR^+$. 
Then $\la _0 \inv \la = u_0 u\inv = 1$.
The map
$$
\pi : G \to \uni \quad \hbox{ given by } \quad \pi (g) = u  \quad (g \in G) 
$$
is a fibration with very rich geometric properties (see
\cite{[PT2]}, \cite{[At]} and the references therein).  We are interested 
in the way that the fibres $\pi \inv (u) = \Gpo u = u \Gpo$ intersect the base
space of a similar fibration induced by a different involution. More precisely, each 
$a \in \Gpo$ induces a $C^*$ involution on $\ca$, namely 
\beq \label{##}
x^{\#_a} = a \inv x^* a.
\eeq
If $\ca$ is represented in the Hilbert space $\cH$, then $a\in \Gpo $ induces
the inner product $\langle , \rangle_a $ given by
$$
\langle \xi , \eta \rangle \sb a = \langle a\xi , \eta \rangle \ , 
\quad \xi , \ \eta \in \cH . 
$$
It is clear that 
$\langle x \xi , \eta \rangle _a = \langle  \xi ,x ^\*a \eta 
\rangle $ for all $x \in \ca$ and $\xi , \eta $ in $\cH$.   
$\ca$ is a \csta \ with this involution and with the norm $\| \cdot \|_a$ 
associated to $\langle , \rangle_a$, $\|x\|_a = \|a\rai x a \mrai \|$, $x \in \ca$. 
For each $a\in \Gpo$, consider the unitary group $\unia$ corresponding to the
involution $\#_a$: 
$$
\unia = \{ g \in G :  g^{\#_a} = g\inv  \} =\{ g \in G : a\inv g^* a = g \inv \}.
$$
We shall study the restriction $\pi \Big{|} _{\unia}$ and the way that 
different $\unia , \ \unib $ are set in $G$. Moreover, we shall also consider 
the \bf $a$-hermitian \rm part of $G$, 
$$
\Gs_a = \{ g\in G : g^{\#_a} = g  \} =\{ g \in G : a\inv g^* a = g  \},
$$
the \bf $a$-positive \rm  part of $\Gs_a $
$$
\Gpo _a = \{ g \in \Gs_a : \sigma (g) \inc \zR^+ \}
$$
and the intersections of these sets when $a$ varies in $\Gpo$. The reader is 
referred to \cite {[L]} and \cite{[D]} for a discussion of operators which 
are hermitian for some inner product.

Observe that each $a\in \Gpo $ induces a fibration
$$
\pi_a : G \to \unia
$$
with fibers homeomorphic to $\Gpo_a$. This paper can be seen 
in some sense as a simultaneous
study of the fibrations $\pi_a$, $a\in \Gpo$.

Let us mention that, from an intrinsic viewpoint, $\unia$ 
can be identified with $
\uni$. Indeed, consider the map $ \fia : \ca \to \ca $ given by 
\beq\label{fia}
\fia (b) = a^{-1/2} b a^{1/2} \quad ( b \in \ca ).
\eeq
Then $\fia (\uni ) = \unia $, $\fia (\Gs ) = \Gs_a$ and 
$\fia (\Gpo ) = \Gpo _a$, since $ \fia : (\ca , * ) \to (\ca , \#_a)$
is an isomorphism of \csta s. We are concerned with the way in 
which the base space and fibers of different fibrations behave 
with respect to each other.

\section{The polar decomposition}

In \cite{[PT1]}, Pedersen and Takesaki proved a technical result which 
was relevant for their generalization of the Sakai's Radon-Nikodym theorem for 
von Neumann algebras \cite{[PT2]}. More precisely, they determined the 
uniqueness and existence of positive solutions of the equation
$$
THT= K
$$
for $H, K $ positive bounded operators in a Hilbert space. We need a weak version
of their result, namely when $H, K$ are positive invertible operators. In this case 
it is possible to give an explicit solution.
\medskip

\begin{lem}[ \cite{[PT1]}]\label{PT1}
If $H, K $ are positive invertible bounded operators in a Hil\-bert space, the equation 
\beq \label{THT}
THT= K
\eeq
has a unique solution, namely 
\beq \label{*}
T = H\mrai ( H \rai K  H \rai ) \rai H \mrai .
\eeq
\end{lem}
\dem
Multiply (\ref{THT}) at left and right by $H\rai$ and factorize 
$$
H\rai THTH\rai = (H\rai T H\rai )^2 .
$$
Then we get the equation 
\beq \label{**}
(H\rai T H\rai )^2 = H\rai K H\rai  .
\eeq
Taking (positive ) square roots and using the invertibility of $H\rai$ we get
the result \QED

\medskip
\noi 
Returning to the map $\pi : G \to \uni$, consider the fiber 
$\pi \inv (u) = \{\la u : \la \in \Gpo \}$. In order to compare the fibration
$\pi$ with $\pi_a$, the following is the key result

\begin{teo} Let $a\in \Gpo$. Then, for every $u \in \uni$ the fiber 
$\pi\inv (u)$ intersects $\unia$ at a single point, namely 
$$
a\mrai ( a \rai u a u\inv a\rai ) \rai a\mrai \cdot u .$$
In other words, the restriction  
$$
\pi | _{\unia} : \unia \to \uni 
$$
is a homeomorphism.
\end{teo}
\dem If $g = \la u \in \unia$ then $ a\inv g^* a = g\inv $ is equivalent
to 
$$
a\inv u\inv \la a = u\inv \la \inv , 
$$
so, after a few manipulations, 
\beq\label{laala}
\la a \la = u a u\inv . 
\eeq
By Pedersen and Takesaki's result, 
there is a unique $\la \in \Gpo $ which satisfies equation
(\ref{laala}) for  fixed $a\in \Gpo$, $u \in \uni$, namely 
\beq\label{pt}
\la = a\mrai ( a \rai u a u\inv a\rai ) \rai a\mrai .
\eeq
Thus, $ (\pi|_{\unia} )\inv : \uni \to \unia $ is given by 
\beq \label{alfa}
(\pi|_{\unia} )\inv (u) = 
a\mrai ( a \rai u a u\inv a\rai ) \rai a\mrai \cdot u
\eeq
which obviously is a continuous map \QED

\medskip
\noi Let $a\in \Gpo$ and consider the involution $\#_a$ defined 
in equation (\ref{##}).
It is natural to look at  those reflections $\eps \in Q$ which 
are $\#_a$-orthogonal, i.e the so called 
$\#_a$-Grassmann manifold of $\ca$. Let us denote this space by
$$
P_a = \{ \eps \in Q : \eps = \eps ^{\#_a} = \eps \inv \} = Q\cap \unia =Q\cap \Gs_a .
$$
In \cite{[PW1]}, Pasternak-Winiarski studied the behavior 
of the orthogonal projection onto a closed subspace of 
a Hilbert space when the inner product varies continuously.
Note that we can identify naturally 
the space of idempotents $q$ with the reflections of $Q$ via 
the affine map $q \mapsto \eps = 2q-1$, which also maps 
the space of orthogonal projections onto $P$.
Based on \cite{[PW1]}, a geometrical study of the space $Q$ 
is made in \cite{[ACS2]}, where the characterization 
$\pi (P_a) = P$ is given 
(proposition 5.1 of \cite{[ACS2]}). In the following proposition we shall 
give a new proof of this fact by showing that the homeomorphism
$\pi|_{\unia} : \unia \to \uni $ maps $P_a \inc \unia $ onto $P \inc \uni$.
Therefore the formula given in equation (\ref{alfa}) for the inverse of $\pi|_{\unia}$
extends the formula given in 
proposition 5.1 of \cite{[ACS2]} for $(\pi|_{P_a})\inv$, 
since they must coincide on $P$.
\medskip
\begin{pro} Let $a\in \Gpo $. Then 
$$\pi (P_a) = \pi ( Q \cap \unia  ) = P .
$$
Therefore $\pi |_{P_a} : P_a \to P$ is a homeomorphism.
\end{pro}
\dem
By the previous remarks, we just need to show that $\pi (P_a ) = P$. 
Observe that if $\eps \in Q$ then $\rho = \pi (\eps ) \in P$: 
in fact, if $\eps = \la \rho $ then $\eps = \eps \inv = \rho \inv 
\la \inv$; but, since the unitary part of $\eps $ corresponding to both
right and left polar decompositions coincide, we get $\rho \inv = 
\rho $. Then $\rho ^* = \rho \inv = \rho $ and $\rho \in P$. Thus, 
$\pi(\unia \cap Q ) \inc P$. 

Let $\alpha = (\pi|_{\unia} )\inv :  \uni \to \unia$. Then by (\ref{alfa}) 
$$
\alpha (u) = a\mrai ( a \rai u a u\inv a\rai ) \rai a\mrai \cdot u.
$$
In order to 
prove the result we need to show that if $\rho \in P$ then $\alpha (\rho ) 
\in P_a = Q\cap \unia $, i.e. $\alpha (\rho ) \in Q$. Indeed, 
%
$$
\barr {rl} 
 \alpha (\rho )^2 & = a\mrai ( a \rai \rho a \rho a\rai ) \rai a\mrai \rho 
a\mrai ( a \rai \rho a \rho a\rai ) \rai a\mrai \rho \\&\\
 & = a\mrai ( (a \rai \rho a\rai)^2 )\rai (a\rai \rho a\rai )\inv 
     ( (a \rai \rho a\rai)^2 ) \rai a\mrai \rho . \earr 
$$
Thus, applying the continuous functional calculus (see e.g. \cite{[Pe]})
to the selfadjoint element $ a \rai \rho a\rai $, if 
$f(t) = |t| = (t^2)\rai $ and $g(t) = {1\over t}$, $t\in \zR \setminus \{ 0\}$, 
$$
\barr{rl}
(\la \rho )^2 & = 
a\mrai f ( a \rai \rho a\rai ) g ( a \rai \rho a\rai ) 
f ( a \rai \rho a\rai ) a\mrai \rho \\&\\
 & = a\mrai [f ( a \rai \rho a\rai )]^2 g ( a \rai \rho a\rai )
a\mrai \rho \\&\\
 & = a\mrai  ( a \rai \rho a\rai ) a\mrai \rho \\&\\
 & = 1 \quad \QED \earr
$$

\bigskip
\begin{num}[Positive parts] \rm
In order to complete the results on the relationship 
between polar decomposition and inner products, 
consider the complementary map of the decomposition
$g = \la u $, namely
\beq\label{pim}
\pim : G \to \Gpo \ , \quad \pim (g) = (gg^*)\rai \ , \quad (g \in G).
\eeq
Of course, there is another ``complementary map'', namely 
$g \mapsto (g^*g)\rai $ corresponding to the decomposition 
$g = u \la '$. We shall see that for every $a \in \Gpo$,  the restriction  
$$
\pim |_{\Gpo _a} : \Gpo _a \to \Gpo 
$$
is a homeomorphism.
Indeed, given $\mu \in \Gpo$, consider the polar
decomposition $ a \mu = \la u$, with $\la \in \Gpo $ and $u \in \uni$. 
Then $ a\inv \la = \mu u^*$, so $\pim ( a\inv \la ) = \mu $ and 
$a \inv \la \in \Gpo _a$, since $a\inv (a\inv \la )^* a = a\inv \la $ and 
the spectrum $\sigma ( a \inv \la ) \inc \zR ^+$. Note that 
$$
(\pim |_{\Gpo _a})\inv (\mu ) = a\inv  \pim (a \mu ),
$$
which is clearly a continuous map. 
An interesting rewriting of the above statement is: 
\end{num}
\bigskip
\begin{pro} If $\ca$ is  a unital \csta \ and 
$a, \la \in \Gpo$, then there exists a unique
$u \in \uni$ such that $$a\la u \in \Gpo.$$
\end{pro}
\dem Indeed, if $g = (\pim |_{\Gpo _a})\inv (\la ) \in \Gpo_a $
and $u = \pi (g)$, then  $\la u =  g \in  \Gpo_a $ means exactly that  
$a\la u \in \Gpo$ \QED

\medskip
It is worth mentioning that  $x \in G$ is the unique positive 
solution of Pedersen-Takesaki equation $xax= b$ if and only if 
$a\rai x b\mrai  \in \uni$. Changing $a, b$ by $a^2 , b^{-2}$ 
respectively, we can write Lemma \ref{PT1} as follows:

\begin{pro} If $\ca$ is  a unital \csta \ and 
$a, b \in \Gpo$, then there exists a unique
$x \in \Gpo$ such that $$axb \in \uni .$$
\end{pro}
\bigskip
\begin{num}[Products of positive operators] \rm
The map $\pou : \Gpo \times \Gpo \to \uni $ given by
$$
\pou(a,b) = axb = a(a\inv(ab^2 a)\rai a\inv)b =(ab^2 a)\rai a\inv b, 
\quad a,b \in \Gpo , 
$$
is not surjective: in fact, the image of $\pou$ consist of those unitary
elements which can be factorized as product of three positive elements. 
On one side $\pou (a, b) = axb \in \uni$ is the product of three 
elements of $\Gpo$. On the other side, if $axb \in \uni$ then by 
Pedersen-Takesaki's result $x$ is the unique positive solution of
$xa^2x = b^{-2}$.

It is easy to show that $-1 \in \uni$ can not 
be decomposed as a product of four positive elements. 
See \cite{[Ph]} and \cite{[W]} for a complete bibliography on these 
factorization problems. See \cite{[Ball]} for more results on 
factorization of elements of $G$ and characterizations of
$P_n = \{ a_1 \dots a_n : a_i \in \Gpo \} $, at least in the finite 
dimensional case.
\end{num}
\bigskip
\begin{num}[Parametrization of the solutions of
Pedersen-Takesaki equations] \rm Given $a,b \in \Gpo$, denote by 
$m = | b\rai a\rai | = (a\rai b a\rai)\rai$. Then the 
set of all solutions of the equation
$xax = b$ is 
$$
\{ a\mrai \ m  \  \eps  \ a\mrai \ : \ \eps \in Q \quad \hbox { and } \quad  
\eps m = m \eps \}.
$$
In fact, $xax=b$ if and only if $(a\rai x a\rai)^2 = m ^2 $ 
and the set of all solutions of $x^2 = c^2 $ for $c\in \Gpo$ is 
$$
\{ c \eps \ : \ \eps ^2 = 1 \quad \hbox { and } \quad \eps c = c \eps \} .
$$
The singular case, which is much more interesting, deserves a particular 
study that we intend to do in several forthcoming papers.
\end{num}



\medskip
\section{Intersections and unions}

\noi For any selfadjoint $c\in \ca$  we shall consider the relative commutant 
sub\csta  
$$
{\cal A}_c = \ca \cap \{ c\} '= \{ d\in {\cal A} : dc=cd \}
$$
and denote by $\uni ({\cal A}_c ) = {\cal A}_c\cap \uni$,  the unitary group of 
$\ca _c$ and, analogously $\Gs (\ca _c)$, $\Gpo (\ca _c)$, $Q(\ca _c)$ and $P(\ca_c )$. 

The space $\Gs$ has a deep relationship with $Q$ (in \cite {[C]} there is a partial 
description of it). Here we only need to notice that the unitary part of any 
$c\in \Gs$ also belongs to $P$. Indeed, if $\la \rho$ is 
the polar decomposition of 
$c$, then $\la \rho = c = c^* = \rho^*\la $. 
By the uniqueness of the unitary part, 
$\rho = \rho^* = \rho\inv \in P$. 
Observe also that $\rho \la = \la \rho$. Moreover, since 
$\la = |c| = (c^2)\inv $, then $\rho = f(c)$ where $f(t) = t \ |t|$. 
So $\rho c = c \rho$.
\medskip
\begin{teo}\label{Gs}
Let $\ca $ a unital \csta \  and $a\in \Gpo$. Then 
$$ \unia \cap \Gs = P(\ca _a ) = \{ u \in P :  au = ua \} .
$$
\end{teo}
\dem By the previous remarks, if $b \in \Gs \cap \unia $  
and $b = \la \rho$ is its polar decomposition, then $\rho \in P$ and
$\rho  \la \inv  = a\inv \rho  \la a $. Using that $\rho \la = \la \rho$
we get  easily 
$$
\la\inv   a \la\inv  = \rho  a \rho  = \la a \la .
$$
By the uniqueness of the positive solution, $\la = \la \inv $ and, since 
$\la \in \Gpo$, this means that $\la = 1 $. Thus $a= \rho a\rho$ and then 
$\rho \in P(\ca _a)$.
Conversely, if $\rho \in 
P(\ca _a )$, then $\rho \in \unia$, since $a\inv \rho ^* a = a\inv \rho a 
= \rho = \rho \inv $ \QED

\begin{rem}\label{varios}\rm 
Let $a \in \Gpo$. Then easy computations show that
\ben
\item \ $\unia \cap \uni  = \uni \cap \ca _a = \uni (\ca _a )$.
\item \ $\Gs _a \cap \Gs  = \Gs \cap \ca _a = \Gs (\ca _a )$.
\item \ $\Gpo _a \cap \Gpo  = \Gpo \cap \ca _a = \Gpo (\ca _a )$.
\item \ $\unia \cap \Gpo =\{1 \}$.
\een
We shall give two proofs of item 4:
 
\noi First proof: $ \pi (\Gpo ) = \{ 1 \} $ but $\pi$ restricted to
$\unia$ is one to one.

\noi Second proof: if $x\in  \unia \cap \Gpo $, then its spectrum
$\sigma (x) \inc \TT \cap \zR^+ = \{1 \}$; on the other side, 
$x$ is normal with respect to the involutions $\#_a$, so $x$ is a normal 
element such that $\sigma (x) = \{ 1 \}$ and it must be $x = 1$.
\end{rem}

\medskip
\noi 
Let $b \in \Gpo$. Recall that the map $\fib $ defined in (\ref{fia})
changes the usual involution by $\#_b$ and also all the corresponding spaces
(e.g $\fib (\Gs ) = \Gs_b$). 

\medskip
\begin{lem}\label{fib} 
Let $b \in \Gpo$. Then for any $c\in \Gpo$, if 
$ d = b\mrai c b\mrai $, 
\beq \label {fibr}
\fib (\uni _d ) = \uni _{c}, \quad 
\fib (\Gs _d )  = \Gs _{c},  \quad \hbox{ and} 
\quad \fib (\Gpo _d )  = \Gpo _{c} 
\eeq
\end{lem}
\dem
Notice that $\uni _d = \fii _d (\uni )$, so $\fib (\uni _d ) = 
\fib \circ \fii _d (\uni )$. But $\fib \circ \fii _d = 
\fii _c \circ Ad (u^* ) $ where $u = \pi (d \rai b \rai )$, 
since $c = |d \rai b \rai |^2 $ and 
$d \rai b \rai =  u c\rai $. As $Ad_{u^*} (\uni ) = \uni$ ( and the same 
happens for $\Gs $ and $\Gpo$),  
we get  $\fib (\uni _d ) = \uni _{c}$ and the other two identities \QED

\medskip
\noi
Then we can generalize the results above for any pair
$a, b \in \Gpo$ instead of $a $ and $ 1$: 

\medskip
\begin{cor} Let $a, b \in \Gpo$ and $c = b\mrai a b\mrai$. Then
\ben
\item \ $\unia \cap \Gs_b =  \fib(\uni _c \cap \Gs ) = 
b\mrai (P (\ca _c ) )b\rai $. 
\item \ $\unia \cap \unib = \fib(\uni _c \cap \uni ) = 
b\mrai \uni (\ca _c ) b\rai.$
\item \ $\Gs _a \cap \Gs _b  = \fib (\Gs_c \cap \Gs) = 
b\mrai\Gs (\ca _c ) b\rai $.
\item \ $\Gpo _a \cap \Gpo _b = \fib( \Gpo_c \cap \Gpo ) = 
b\mrai \Gpo (\ca _a ) b\rai $.
\item \ $\unia \cap \Gpo _b = \fib (\uni_c \cap \Gpo ) = \{1 \}$.
\een 
\end{cor}
\dem Use Proposition \ref{Gs}, Remark \ref{varios} and Lemma \ref{fib}
\QED

\bigskip

\noi In the following proposition we describe the set of 
elements of $\ca$ which are unitary (resp. positive, Hermitian) for some
involution $*_a$  $(a \in \Gpo )$. We state the result without proof.

\begin{pro} If $\ca$ is a unital \csta , the following identities hold: 
\beq 
\bigcup_{a \in \Gpo} \unia 
= \bigcup_{g\in G} g\ \uni g\inv 
= \bigcup_{a\in \Gpo} a \ \uni a\inv , 
\eeq
\beq
\bigcup_{a \in \Gpo} \Gpo _a 
= \bigcup_{g\in G} g\ \Gpo  g\inv 
= \bigcup_{a\in \Gpo} a \ \Gpo  a\inv = \Gpo \Gpo ,
\eeq
where $\Gpo \Gpo  = \{ ab : a,b\in \Gpo \}$ and
\beq
\bigcup_{a \in \Gpo} \Gs _a 
= \bigcup_{g\in G} g\ \Gs  g\inv 
= \bigcup_{a\in \Gpo} a \ \Gs  a\inv = \Gs \Gpo = \Gpo \Gs.
\eeq
\end{pro}
The following example shows that there is no obvious 
spectral characterization of these subsets of $G$: if $x$ is
nilpotent, then $1+x$ does not belong to any of them
but $\sigma (1+x) = \{1\} \inc \zR , \zR^+ , \TT$.

\begin{num}[Final geometrical remarks] \rm
The subsets of $\ca $ studied in this paper have all a rich structure
as differential manifolds. The reader is referred to 
\cite{[Ha]} and \cite{[At]} for the case of $\uni$ and to \cite{[C]} (and
the references therein) for $Q, P, \Gs $ and $\Gpo$. The map $\fia$ defined
in equation (\ref{fia}) is 
clearly a diffeomorphism which allows to get all the information on 
$\unia , \Gs _a , \Gpo _a $ from that available on $\uni , \Gs , \Gpo$, 
respectively. The main results of the paper say that the map 
$\pi$ is a diffeomorphism between $\unia $ and $\uni$, $P_a $ and $P$
and so on.
\end{num}

\vglue1truecm

\end{document}